\documentclass[11pt,reqno]{amsart}
\usepackage{graphicx}
\usepackage{graphics}
\usepackage{amsmath}
\usepackage{epsfig,bm}
\usepackage{rotating}
\usepackage{url}
\usepackage{float} 
\usepackage{amsthm}

\usepackage{enumitem}

\DeclareMathOperator{\perm}{perm}
\DeclareMathOperator{\LS}{LS}
\DeclareMathOperator{\RS}{RS}

\begin{document}
\title[Expectation Value of Permanent Products]{Asymptotic Behavior of the Expectation Value of Permanent Products, a Sequel}

\author[P. Federbush]{P. Federbush\\
Department of Mathematics\\
	University of Michigan \\
	Ann Arbor, MI 48109-1043, USA}

\maketitle
\begin{abstract}

Continuing the computations of the previous paper, \cite{1}, we calculate another approximation to the expectation value of the product of two permanents in the ensemble of 0-1 $n\times n$ matrices with like row and column sums equal $r$ uniformly weighted. Here we consider the Bernoulli random matrix ensemble where each entry independently has a probability $p=r/n$ of being one, otherwise zero. We denote the expectations of the approximation ensemble of \cite{1} by $E$, and the expectations of the present approximation ensemble, the Bernoulli random matrix ensemble, by $\bar{E}$. One has for these
$$
	\lim_{r\to\infty}\left(\lim_{n\to\infty}\frac{1}{n}
	\ln(E(\perm_m(A)))-\lim_{n\to\infty}\frac{1}{n}
	\ln(\bar{E}(\perm_m(A)))\right)=0
$$
and
\begin{multline*}
	\lim_{n\to\infty}\frac{1}{n}\ln(E(\perm_m(A)\perm_{m'}(A)))=\\
	\lim_{n\to\infty}\frac{1}{n}\ln(E(\perm_m(A)))+\lim_{n\to\infty}\frac{1}{n}
	\ln(E(\perm_{m'}(A)))
\end{multline*}
Here and in all such formulas the subscripts $m,m'$ are assumed proportional to $n$. It seems likely to us that
\begin{multline*}
	\lim_{r\to\infty}\left(\lim_{n\to\infty}\frac1n\ln(\bar{E}(\perm_m(A)\perm_{m'}(A)))\right.\\
	\left.-\lim_{n\to\infty}\frac1n\ln(\bar{E}(\perm_m(A)))-\lim_{n\to\infty}\frac1n
	\ln(\bar{E}(\perm_{m'}(A)))\right)=0
\end{multline*}
We believe: ``$E$ gives us the `correct' expectations in these equations, and $\bar{E}$ is only `correct' in the $r\to\infty$ limit.''
\end{abstract}

\newpage

In the previous paper in this series, \cite{1}, we worked with a measure on matrices giving uniform weight to a sum of $r$ independent 
random $n\times n$ permutation matrices and computed
\begin{multline}\label{1}
	\lim_{n\to\infty}\frac1n\ln(E(\perm_m(A)\perm_{m'}(A)))=\\
	\lim_{n\to\infty}\frac1n\ln(E(\perm_m(A)))+\lim_{n\to\infty}\frac1n\ln(E(\perm_{m'}(A)))
\end{multline}
where $r$ is fixed, and $m$ and $m'$ are each proportional to $n$, the size of the matrices. We expect this to be the same result as obtained using the `correct' measure, the uniform measure on 0-1 matrices with row and column sums all equal $r$. This is the result suggested by Friedland's Asymptotic Matching Conjecture, now proven, \cite{2}, \cite{3}, \cite{4}. The computation of this note is in the measure on all 0-1 matrices where the entries are independently all zeros and ones, with the probability of being one being $p=r/n$. It does seem suggested by numerical study that \eqref{1} is true in the limit $r\to\infty$ for the current computation, see table at the end of paper. Unlike the computation in \cite{1} it seems easy to make the current computation rigorous. The impossible dream for future work is getting a handle on expectations of arbitrary finite permanent products.

We first compute the expectation of the permanent of a single matrix in the Bernoulli random matrix ensemble, using a bar to distinguish the expectations of this paper from those in \cite{1}.
\begin{equation}\label{2}
	\bar{E}(\perm_m(A))=\binom{n}{m}\binom{n}{m}m!\,p^m
\end{equation}
we let $m=an$, and have $p=r/n$, so 
\begin{equation}\label{3}
	\bar{E}(\perm_m(A))=\binom{n}{an}\binom{n}{an}(an)!\,(r/n)^m
\end{equation}
This differs from the corresponding result in \cite{1}, but looking at the large $r$ limit it is easy to show 
\begin{equation}\label{4}
	\lim_{r\to\infty}\left(\lim_{n\to\infty}\frac1n\ln(\bar{E}(\perm_m(A)))-
	\lim_{n\to\infty}\frac1n\ln(E(\perm_m(A)))\right)=0
\end{equation}
where $m$ is understood to be $an$ also in the second expectation here. We expect a similar equation to hold for expectations of a product of two permanents.

We proceed to study
\begin{equation}\label{5}
	\bar{E}(\perm_m(A)\perm_{m'}(A))
\end{equation}
We use the relationship between $n\times n$ 0-1 matrices and bipartite graphs on $n$ black and $n$ white vertices in describing our computation. The calculation of \eqref{5} is done by counting the number of pairs of an $m$-matching and an $m'$-matching on the complete $n\times n$ bipartite graph, and multiplying each such term by the probability that all the edges in
the two matchings are present. We introduce non-negative integers $A,B,C,\bar A,\bar B,D,Q,Z,W,X,E$ that we will associate to each term (an $m$-matching and an $m'$-matching). They are specified as follows:
\begin{enumerate}[label=\arabic*)]
\item $C$ is the number of common black vertices of the two matchings.

\item $D$ is the number of common white vertices of the two matchings.

\item $A$ is the number of black vertices in the $m$-matching not present in the $m'$-matching, so 
\begin{equation}\label{6}
	A+C=m
\end{equation}

\item $B$ is the number of black vertices in the $m'$-matching not present in the $m$-matching, so
\begin{equation}\label{7}
	B+C=m'
\end{equation}

\item $\bar A$ and $\bar B$ are similarly defined for the white vertices, so
\begin{equation}\label{8}
	\bar A+D=m
\end{equation}
\begin{equation}\label{9}
	\bar B+D=m'
\end{equation}

\item $Q$ is the number of edges shared by the $m$-matching and the $m'$-matching. 
Without loss of generality we may assume 
\begin{equation}\label{10}
	Q\leq C\leq D\leq m\leq m'
\end{equation}
For each of our integers $A,B,\dots,E$, we use the subscript $s$ to indicate the set of vertices (or edges) counted by the integer. For example, $\bar A_s$ is the set of white vertices in the $m$-matching not present in the $m'$-matching.

\item $Z$ is the number of vertices in $C_s$ such that the edge in the $m$-matching with this vertex (in its boundary) and the edge in the $m'$-matching with this vertex are not identical but both have their other vertices in $D_s$.

\item $W$ is the number of vertices in $C_s$ such that the edge in the $m$-matching with this vertex has its other vertex in $D_s$, but the edge in the $m'$-matching with this vertex has its other vertex not in $D_s$.

\item $X$ has the same definition with the roles of $m$ and $m'$ interchanged.

\item $E$ is the number of vertices in $C_s$ such that both edges containing this vertex have other vertices not in $D_s$.
\end{enumerate}

We note 
\begin{equation}\label{11}
Q+Z+W+X+E=C
\end{equation}

We compute \eqref{5} by fixing the parameters $A,B,\dots,E$ and then summing over the number of $m$-matchings and $m'$-matchings consistent with these values, and then summing over the values of the parameters. To compute
\begin{equation}\label{12}
\lim_{n\to\infty}\frac1n\ln(\bar E(\perm_m(A)\perm_{m'}(A)))
\end{equation}
one need only keep the largest term in the sum over the parameters. 

We turn to the sum over $m,m'$-matchings consistent with a fixed choice of parameters. We write this sum as
\begin{equation}\label{13}
\prod_{i=1}^{11}T_i
\end{equation}
and turn to specifying the $T_i$, and where they arise.

\begin{enumerate}[label=$T_\text{\arabic*}$)]
\item $T_1$ specifies the probability that the edges in the $m$-matching and $m'$-matching are all present. Since there are $m+m'-Q$ edges present
\begin{equation}\label{14}
T_1=p^{m+m'-Q}=\left(\frac{r}{n}\right)^{m+m'-Q}
\end{equation}

\item $T_2$ selects the black vertices in $A_s$, $C_s$, and $B_s$
\begin{equation}\label{15}
T_2=\frac{n!}{A!\,C!\,B!\,(n-A-B-C)!}
\end{equation}

\item $T_3$ similarly selects the white vertices in $\bar A_s, D_s,\bar B_s$
\begin{equation}\label{16}
T_3=\frac{n!}{\bar A!\,D!\,\bar B!\,(n-\bar A-D-\bar B)!}
\end{equation}

\item $T_4$ divides the black vertices in $C_s$ into these in $Z_s,W_s,X_s,E_s,$ and the black vertices of the edges in $Q_s$.
\begin{equation}\label{17}
T_4=\frac{C!}{Q!\,Z!\,W!\,X!\,E!}
\end{equation}

\item $T_5$ selects the edges in $Q_s$
\begin{equation}\label{18}
T_5=\frac{D!}{(D-Q)!}
\end{equation}

\item $T_6$ selects the edges whose black vertices are in $Z_s$.
\begin{equation}\label{19}
T_6\cong\frac{(D-Q)!}{(D-Q-Z)!}\frac{(D-Q)!}{(D-Q-Z)!}
\end{equation}
We write the approximate equality since we have not taken into account in the right that the two edges leaving the same black vertex must be distinct. Neglecting this requirement does not matter in the limit in \eqref{12}.

\item In $T_7$ we select the edges whose black vertex lies in $W_s$ and whose white vertex lies in $D_s$
\begin{equation}\label{20}
T_7=\frac{(D-Q-Z)!}{(D-Q-Z-W)!}
\end{equation}

\item In $T_8$ we select the edges whose black vertex lies in $X_s$ and whose white vertex lies in $D_s$
\begin{equation}\label{21}
T_8=\frac{(D-Q-Z)!}{(D-Q-Z-X)!}
\end{equation}
One may find it difficult to convince oneself that the enumerations of \eqref{20} and \eqref{21} are correct.

\item In $T_9$ we select the edges whose black vertices lie in $E_s$ or $W_s$ and whose white vertices lie in $\bar{B}_s$, as well as the edges whose black vertices lie in $E_s$ or $X_s$ and whose white vertices in $\bar A_s$
\begin{equation}\label{22}
T_9=\frac{\bar B!}{(\bar B-W-E)!}\cdot\frac{\bar A!}{(\bar A-X-E)!}
\end{equation}

\item In $T_{10}$ we select the edges whose black vertices lie in $A_s$
\begin{equation}\label{23}
T_{10}=A!
\end{equation}

\item
In $T_{11}$ we select the edges whose black vertices lie in $B_s$
\begin{equation}\label{24}
T_{11}=B!
\end{equation}

\end{enumerate}

We turn to studying 
\begin{equation}\label{25}
\lim_{n\to\infty}\frac1n\ln(\bar E(\perm_m(A)\perm_{m'}(A)))
\end{equation}
where
\begin{equation}\label{26}
m=sn
\end{equation}
and
\begin{equation}\label{27}
m'=tn
\end{equation}
we also scale $A,B,\dots,E$ using lower case letters to have 
\begin{equation}\label{28}
A=na,B=nb,\cdots,E=ne
\end{equation}
Equations \eqref{6}--\eqref{11} become
\begin{equation}\label{29}
a+c=s
\end{equation}
\begin{equation}\label{30}
b+c=t
\end{equation}
\begin{equation}\label{31}
\bar a+d=s
\end{equation}
\begin{equation}\label{32}
\bar b+d=t
\end{equation}
\begin{equation}\label{33}
q\leq c\leq d\leq s\leq t
\end{equation}
\begin{equation}\label{34}
q+z+w+x+e=c
\end{equation}

We fix $s$ and $t$ and parameters $a,b,\dots,e$ consistent with \eqref{29}--\eqref{34}. We evaluate 
\begin{equation}\label{35}
\frac1n\ln\!\left(\prod_{i=1}^{11}T_i\right)
\end{equation}
with the $T_i$ given in \eqref{14}--\eqref{24}, and approximating 
\begin{equation}\label{36}
\ln f!\cong f\ln f-f
\end{equation}
Then
\begin{equation}\label{37}
\lim_{n\to\infty}\frac1n\ln(\bar E(\perm_m(A)\perm_{m'}(A)))=\max\!\left(\frac1n\ln\!\left(\prod_{i=1}^{11}T_i\right)\right)
\end{equation}
where the maximum is over all compatible values of the parameters $a,b,\dots,e$. (Equations \eqref{29}--\eqref{34} are imposed.)

We have done a little numerical study with $s=t$. We set
\begin{equation}\label{28}
\LS=\lim_{n\to\infty}\frac1n\ln(\bar E(\perm_{m}(A)\perm_{m}(A)))
\end{equation}
\begin{equation}\label{39}
\RS=2\lim_{n\to\infty}\frac1n\ln(\bar E(\perm_m(A)))
\end{equation}
and computed the following table
\begin{center}
\begin{tabular}{c|c|c}
&$\LS$&$\RS$\\
\hline
$s=.4$, $r=5$&2.4771&2.4466\\
$s=.4$, $r=50$&4.2918&4.2886\\
$s=.8$, $r=5$&2.7435&2.6197\\
$s=.8$, $r=50$&6.3166&6.3038\\
\end{tabular}
\end{center}
To obtain those results we asked Maple to solve numerically the algebraic equations obtained setting the derivatives of \eqref{35} with respect to a full set of free parameters equal to zero. It seems reasonable to conjecture $\LS$ and $\RS$ converge to each other as $r\to\infty$, and even likely that this is not hard to prove.


\end{document}